\documentclass{amsart}
\usepackage{amssymb,latexsym}
\usepackage{stmaryrd}
\usepackage{xcolor}

\newtheorem{theorem}{Theorem}[section]
\newtheorem{lemma}[theorem]{Lemma}

\theoremstyle{definition}

\numberwithin{equation}{section}

\newcommand{\D}{{\mathbb D}}

\newcommand{\T}{{\mathbb T}}

\newcommand{\HM}{\mathcal{M}}

\newcommand{\HH}{\mathcal{H}}

\makeatletter
\@namedef{subjclassname@2020}{%
  \textup{2020} Mathematics Subject Classification}
\makeatother

\begin{document}

\date{}

\author[S. Luo]{Shuaibing Luo}
\email{sluo@hnu.edu.cn}
\address{School of Mathematics, Hunan University, Changsha, 410082, PR China}

\subjclass[]{}

\title[]{A note on the equality between de Branges-Rovnyak and Dirichlet spaces}

\begin{abstract}
In this note, we demonstrate that the methods in \cite{BDL} can be used to characterize when a de Branges-Rovnyak space coincides with a
Dirichlet-type space.

This note is not intended for publication.
\end{abstract}
\maketitle

%\tableofcontents

\section{Introduction}

\subsection{de Branges-Rovnyak spaces}
For $b$ in $H^\infty_1$, the de Branges--Rovnyak space $\mathcal{H}(b)$ is the
reproducing kernel Hilbert space of analytic functions on $\mathbb D$ with kernel
\[
    k_w^b(z)=\frac{1-\overline{b(w)}\,b(z)}{1-\overline wz},\qquad z,w\in\mathbb D .
\]

The structure of $\mathcal{H}(b)$ depends on the properties of the function $b$.
When $b$ is non-extreme there exists a unique outer function \(a\), normalized by \(a(0)>0\), such that
\[
|a|^2+|b|^2=1
\qquad m\text{-almost everywhere on }\mathbb T.
\]
The function \(a\) is called the Pythagorean mate of \(b\), and we refer to $(b,a)$ as a Pythagorean pair. When the Pythagorean pair $(b,a)$ satisfies the corona condition, $\mathcal{H}(b)$ coincides with the range of the co-analytic Toeplitz operator $T_{\bar{a}}$, that is,
\[
    \mathcal{H}(b)=\mathcal M(\bar a):=\{f\in H^2:\ f=T_{\bar a}g\ \text{for some }g\in H^2\},
\]
with equivalent norms, where
\(
    \|T_{\bar a}g\|_{\mathcal M(\bar a)}=\|g\|_{H^2}
\) (the operator $T_{\bar{a}}$ is injective, since $a$ is outer). For more details on de Branges--Rovnyak spaces, we refer to \cite{FM16, Sa94}.

\subsection{Dirichlet-type spaces}
Suppose $\mu$ is a finite positive Borel measure on the closed unit disc $\overline{\D}$. Let
\begin{equation*}
\label{e11}
U_\mu(z) = \int_{\D}\log\left|\frac{1-\overline{w}z}{z-w}\right|^2\frac{d\mu(w)}{1-|w|^2}+\int_{\T}\frac{1-|z|^2}{|1-\overline{w}z|^2}d\mu(w),\quad  z\in \D.
\end{equation*}
Then $U_\mu(z)$ is a positive superharmonic function on $\D$. The Dirichlet-type space $D(\mu)$ consists of holomorphic functions $f$ on $\D$ such that
$$\int_\D |f'(z)|^2 U_\mu(z) dA(z)<\infty.$$
It is konwn that $D(\mu)$ is contained in the Hardy space $H^2$, and the norm of $f$ in $D(\mu)$ is defined by
$$\|f\|^2_{D(\mu)} = \|f\|_{H^2}^2 + \int_\D |f'(z)|^2 U_\mu(z) dA(z).$$
Note that there is a one to one correspondence between the set of finite positive Borel measures $\mu$ on $\overline{\D}$ and the collection of positive superharmonic functions $U_\mu$ on $\D$ (\cite{La72}). For the measures $\mu$ supported on the unit circle $\T$ the Dirichlet-type space was introduced in \cite{Ri91}. The general case that $\mu$ on $\overline{\D}$ was investigated in \cite{Al93}. We refer to \cite{Al93, EKMR14} for the background of the $D(\mu)$ spaces.

For $f \in H^2$ and $z \in \overline{\D}$ the local Dirichlet integral of $f$ at $z$ is defined by
$$D_z(f) = \int_\T \left|\frac{f(w) - f(z)}{w-z}\right|^2 \frac{|dw|}{2\pi}.$$
If $f$ does not have a nontangential limit at $z \in \T$, then $D_z(f)$ is defined to be $\infty$. By \cite[Theorem IV.1.9]{Al93} we have
\begin{align*}
\int_\D |f'(z)|^2 U_\mu(z)\,dA(z)=\int_{\overline{\D}}D_z(f)\,d\mu(z).
\end{align*}
Hence
$$D(\mu) = \left\{f \in H^2: \int_{\overline{\D}}D_z(f)\,d\mu(z) < \infty\right\}.$$

The investigation of the relationship between $\HH(b)$ and $D(\mu)$ has a long history and the relationship has been studied by many authors. The paper \cite{BDL} has an overview of this history. When $\mu$ is supported on the unit circle $\T$, the authors in \cite[Theorem 1.3]{BDL} gave a complete answer to the question of when $\HH(b) = D(\mu)$ holds with equivalent norms. The solution is given in terms of a corona condition and a potential condition.

In this note, we demonstrate that the methods in \cite{BDL} can also be used to characterize when $\HH(b) = D(\mu)$ for $b \in H^\infty_1$ and finite positive Borel measure $\mu$ on the closed unit disk.

Let $\mu$ be a finite positive Borel measure on $\overline{\D}$. Define
\begin{equation*}
V_\mu(z) = \int_{\overline{\D}} \frac{d\mu(\lambda)}{|1 - \overline{\lambda}z|^2}, \qquad z \in \mathbb{D}.
\end{equation*}

\begin{theorem}\label{mainthmduhb}
Let \(b \in H^\infty_1\), and let \(\mu\) be a finite positive Borel measure on \(\overline{\mathbb{D}}\). Then
$
\mathcal{H}(b) = D(\mu)
$
if and only if the following conditions hold:
\begin{enumerate}
    \item \(b\) is non-extreme, and the Pythagorean pair $(b,a)$ satisfies the corona condition, that is
    \begin{equation*}
    \inf_{z \in \mathbb{D}} (|a(z)| + |b(z)|) > 0.
    \end{equation*}

    \item There are constants \(c, C > 0\) such that
    \begin{equation*}
    c \le |a(z)|^2(1 + V_\mu(z)) \le C, \qquad z \in \mathbb{D}.
    \end{equation*}
\end{enumerate}
\end{theorem}

\section{Proof of Theorem \ref{mainthmduhb}}
The proof is similar to the proof in \cite{BDL}. We sketch the proof in the following.

Let \(b \in H^\infty_1\). For $0 < r <1$, let
\[
\Omega_r = \{z \in \mathbb{D} : |b(z)| < r\}
\]
be the sublevel set.

When $\mu$ is supported on $\T$, the following lemma is obtained in \cite[Lemma 3.2]{BDL}. When $\mu$ is supported on $\D$, the calculation is straightforward.

\begin{lemma}\label{sublevelbdvmu}
Let \(b \in H^\infty_1\) (not necessarily non-extreme), and let \(\mu\) be a finite positive Borel measure on \(\overline{\mathbb{D}}\). Suppose that \(\mathcal{H}(b) \subseteq D(\mu)\).
Let $0 < r <1$. If \(\Omega_r\) is nonempty, then
\[
\sup_{w \in \Omega_r} V_\mu(w) < \infty.
\]
\end{lemma}

\begin{proof}
Let $\sigma=\mu|_{\mathbb{T}}$ and $\nu=\mu|_{\mathbb{D}}$. Then $V_\mu = V_\sigma + V_\nu$. Since $\mathcal{H}(b)\subseteq D(\mu)\subseteq D(\sigma)$, \cite[Lemma 3.2]{BDL} gives
\[
\sup_{w \in \Omega_r} V_\sigma(w) < \infty.
\]

Now we show $\sup_{w \in \Omega_r} V_\nu(w) < \infty$.
For $\lambda \in \D$, let
\begin{align*}
Q_\lambda f(z) = \frac{f(z) - f(\lambda)}{z - \lambda}.
\end{align*}
Since $\mathcal{H}(b)\subseteq D(\mu)\subseteq D(\nu)$,
there exists a constant $C>0$ such that
\begin{align}\label{embeddingestkel}
\int_{\mathbb{D}}\|Q_\lambda f\|_{H^2}^2\,d\nu(\lambda)\leq C\|f\|_{\mathcal{H}(b)}^2,\qquad f\in\mathcal{H}(b).
\end{align}

Take $w\in\Omega_r$ with $|w|\geq(1+r)/2$ and put $\beta=b(w)$. We have $(1-\bar{w}z)K_w^b(z)=1-\bar{\beta}b(z)$.
Then for $\lambda\in\mathbb{D}$, by calculation, we obtain that
\[
(1-|w|^2)(Q_\lambda K_w^b)(w)=\bar{w}K_w^b(\lambda)-\bar{\beta}(Q_\lambda b)(w).
\]
By Schwarz-Pick inequality, we have
\begin{align*}
|(Q_\lambda b)(w)| = \frac{|b(w) - b(\lambda)|}{|w-\lambda|} \leq \frac{|1-\overline{b(\lambda)}b(w)|}{|1-\overline{\lambda}w|} = |K_w^b(\lambda)|.
\end{align*}
Hence
\[
(1-|w|^2)|(Q_\lambda K_w^b)(w)|\geq\big(|w|-|\beta|\big)|K_w^b(\lambda)|\geq\frac{\big(|w|-|\beta|\big)(1-|\beta|)}{|1-\bar{\lambda}w|}.
\]
So
\begin{align*}
\frac{\big(|w|-|\beta|\big)(1-|\beta|)}{|1-\bar{\lambda}w|}& \leq (1-|w|^2)|(Q_\lambda K_w^b)(w)|\\
& \leq (1-|w|^2) \|Q_\lambda K_w^b\|_{H^2} \frac{1}{\sqrt{1-|w|^2}}\\
& = \sqrt{1-|w|^2} \|Q_\lambda K_w^b\|_{H^2}.
\end{align*}
It follows that
\begin{align*}
\frac{\|Q_\lambda K_w^b\|_{H^2}^2}{\|K_w^b\|_{\mathcal{H}(b)}^2} &\geq \frac{\big(|w|-|\beta|\big)^2(1-|\beta|)^2}{|1-\bar{\lambda}w|^2}\frac{1}{(1-|w|^2)\|K_w^b\|_{\mathcal{H}(b)}^2}\\
& = \frac{\big(|w|-|\beta|\big)^2(1-|\beta|)^2}{|1-\bar{\lambda}w|^2} \frac{1}{1-|\beta|^2}\\
& = \frac{\big(|w|-|\beta|\big)^2(1-|\beta|)}{|1-\bar{\lambda}w|^2(1+|\beta|)}\\
& \geq\frac{(1-r)^3}{8|1-\bar{\lambda}w|^2}.
\end{align*}
Thus by (\ref{embeddingestkel}), we conclude that $V_\nu(w)\leq 8C/(1-r)^3$.
\end{proof}

When $\mu$ is supported on $\T$, the following result is obtained in \cite[Theorem 3.3]{CR13}. We use the same argument here.
\begin{lemma}\label{estimate1plusvmu}
Let \(b \in H^\infty_1\), and let \(\mu\) be a finite positive Borel measure on \(\overline{\mathbb{D}}\).
Suppose that \(\mathcal{H}(b) = D(\mu)\), and $(b,a)$ is the Pythagorean pair.
Then
\begin{align*}
1 + \frac{|b(w)|^2}{|a(w)|^2} \asymp 1 + V_\mu(w), \qquad w \in \D.
\end{align*}
\end{lemma}
\begin{proof}
Let \(k_w(z) = (1 - \overline{w}z)^{-1}\) be the Hardy kernel. Then
\[
Q_\lambda k_w = \frac{\overline{w}}{1 - \overline{w}\lambda} k_w.
\]
So
\begin{equation*}
\|k_w\|_{D(\mu)}^2 = \|k_w\|_{H^2}^2 + \int_{\overline{\D}} \|Q_\lambda k_w\|_{H^2}^2 d\mu(\lambda) = \frac{1 + |w|^2 V_\mu(w)}{1 - |w|^2}.
\end{equation*}
\cite[Lemma 3.4]{CR13} gives
\begin{equation*}
\|k_w\|_{\mathcal{H}(b)}^2 = \frac{1 + |b(w)/a(w)|^2}{1 - |w|^2}.
\end{equation*}
Thus the conclusion follows from $\|k_w\|_{\mathcal{H}(b)}^2 \asymp \|k_w\|_{D(\mu)}^2, w \in \D$.
\end{proof}

When $\mu$ is supported on $\T$, the following result is obtained in \cite[Theorem 1.2]{BDL}. We use the same argument here.
\begin{lemma}\label{nececondcorona}
Let \(b \in H^\infty_1\), and let \(\mu\) be a finite positive Borel measure on \(\overline{\mathbb{D}}\).
Suppose that \(\mathcal{H}(b) = D(\mu)\), and $(b,a)$ is the Pythagorean pair. Then $(b,a)$ satisfies the corona condition.
\end{lemma}
\begin{proof}
Let $0 < r < \sqrt{2}/2$. Lemma \ref{sublevelbdvmu} implies that $V_\mu$ is bounded on the sublevel set $\Omega_r$. Then Lemma \ref{estimate1plusvmu} ensures that
\(b/a\) is bounded on $\Omega_r$. Thus by \cite[Lemma 3.1]{BDL}, we obtain that $(b,a)$ satisfies the corona condition.
\end{proof}

Now we can prove Theorem \ref{mainthmduhb}.
\begin{proof}[Proof of Theorem \ref{mainthmduhb}]
The necessity follows from Lemmas \ref{estimate1plusvmu} and \ref{nececondcorona}.

Note that if $(b,a)$ satisfies the corona condition, then $\HH(b) = \HM(\overline{a})$ (\cite{FM16, Sa94}).
Thus to prove the sufficiency part, it suffices to prove that $D(\mu) = \HM(\overline{a})$. This follows from the proof of \cite[Theorem 1.3]{BDL}. \cite[Theorem 1.3]{BDL} deals with measures $\mu$ supported on $\T$. But from the proof of \cite[Theorem 1.3]{BDL}, we see that the same proof works for measures supported on $\overline{\D}$.
\end{proof}

%\noindent \textbf{Acknowledgement}:
%This work was supported by NSFC (No.12271149), Natural Science Foundation of Hunan Province (No.2024JJ2008).

\end{document}